\documentclass[12pt]{article}
\usepackage[cp1251]{inputenc}
\usepackage[russian]{babel}
\usepackage{amssymb,amsmath,amsfonts,latexsym,mathtext}
\begin{document}
\author{Белошапка В.К.}

\date{}

\title{{\bf О группе голоморфных автоморфизмов модельной поверхности}}

\maketitle

\begin{abstract}
 В работе доказано, что группа ${\rm Aut}\, Q$  всех голоморфных автоморфизмов голоморфно однородной невырожденной  модельной поверхности $Q$ представляет собой подгруппу группы бирациональных изоморфизмов объемлющего комплексного пространства (группы Кремоны) равномерно ограниченной степени. Дается оценка степени через размерность объемлющего пространства (теорема 4). Показано, что ни одно из условий теоремы нельзя  ослабить. Также в работе рассмотрен вопрос о связности ${\rm Aut}\, Q$ (теорема 7). Данная работа непосредственно примыкает к предыдущей работе автора \cite{VB20}.
\end{abstract}

\footnote{
Механико-математический факультет Московского университета им.Ломоносова,
Воробьевы горы, 119992 Москва, Россия, vkb@strogino.ru.

Работа выполнена  при финансовой   поддержке гранта РНФ 18-41-05003}

{\bf Введение.} В работах \cite{VB04} и \cite{VB20}, а также других работах автора $CR$-многообразия изучались с помощью единообразного подхода - {\it метода модельной поверхности}. При этом класс рассматриваемых многообразий определялся с помощью некоторого условия "невырожденности". В \cite{VB04} - это было условие полной невырожденности ростка $CR$-многообразия. Условие полной невырожденности является условием общего положения в том смысле, что любой росток положительной $CR$-размерности малой гладкой деформацией можно сделать вполне невырожденным. Однако тип такого ростка по Блуму-Грэму не может быть произвольным. В работе \cite{VB20} та же программа метода модельной поверхности (стандартный набор утверждений) была реализована для класса $CR$-ростков произвольного {\it конечного} Блум-Грэм-типа с условием {\it голоморфной невырожденности}. При выполнении этой пары условий мы называем росток {\it невырожденным}. В упомянутый стандартный набор утверждений входит утверждение о бирациональности любого автоморфизма модельной поверхности. В работе \cite{VB20} это утверждение было пропущено.
Оно было анонсировано в списке открытых вопросов (гипотеза 5). В данной работе этот пропуск устранен (теорема 4), т.е. доказано, что бирациональность и ограниченность степени остаются при соблюдении двух условий:  невырожденности и голоморфной однородности. В прежних версиях этого утверждения голоморфная однородность не фигурировала в качестве самостоятельного условия. Это связано с тем, что все вполне невырожденные модельные поверхности являются однородными. В новом контексте невырожденных поверхностей - это уже не так.  В \cite{VB20} был дан критерий однородности модельной поверхности. А именно, было показано, что постоянство Блум-Грэм-типа, которое, очевидно, является необходимым условием однородности - для модельной поверхности является достаточным.  Далее там же было показано, что набор весов однородной модельной поверхности совершенно не произволен. А именно, было показано, что этот набор имеет следующий вид:
$m_1=2,\, m_2=3, \dots, m_l=l+1$.

\vspace{5ex}

Первое утверждение о бирациональности автоморфизмов модельной поверхности было доказано  А.Тумановым \cite{T}
с помощю приема В.Каупа \cite{K}. Это было сделано тогда для простейшего Блум-Грэм-типа, а именно для $m=(2,k)$ (модельная квадрика коразмерности $k$). Наше построение -- это, по существу, многократное рекурсивное использование рассуждения Каупа-Туманова. Отметим при этом, что возможность доказательства бирациональности для модельной поверхности произвольного Блум-Грэм-типа по такой схеме связана с важной структурной особенностью модельных поверхностей -- с их "треугольностью".

\vspace{3ex}

Второй основной результат данной работы (теорема 7) - описание  топологического строения группы голоморфных автоморфизмов модельной поверхности. Этот результат является новым лишь отчасти. Во всех версиях метода модельной поверхности, начиная с \cite{VB91},  подразумевалось следующее простое следствие базовых конструкций. Поскольку подгруппа  $\mathcal{G}_{+}$ (нелинейные автоморфизмы, сохраняющие начало координат) параметризуется ядром гомологического оператора, то она, как всякое линейное пространство, является связной и односвязной (пункт (a) теоремы). Мы приводим здесь доказательство по двум причинам. Первая: начиная с работы \cite{VB20} метод модельной поверхности работает в гораздо более широком контексте. Вторая: наличие явно сформулированного утверждения - это возможность для ссылок.  Остальное содержание теоремы   (пункты (b) и (c)) - являются новыми.

\vspace{5ex}

{\bf  Бирациональность.} Чтобы сделать ясным доказательство общего утверждения, мы предпошлем ему разбор одного частного случая.
А именно, мы рассмотрим первый тип, выходящий за рамки квадратичных моделей $m=\{(2,k),(3,K)\}$. При этом мы не будем предполагать, что поверхность  является вполне невырожденной.   Для вполне невырожденного случая кратность $k$ должна быть равной $n^2$, где $n$ -- это $CR$-размерность. В общем же случае мы имеем $1 \leq k \leq n^2$.  Невырожденные модельные поверхности указанного типа были подробно рассмотрены в \cite{VB18} и результат о бирациональности там был анонсирован.

\vspace{3ex}

Итак,  модельная поверхность $Q$ Блум-Грэм-типа $m=\{(2,k),(3,K)\}$ -- это поверхность в пространстве ${\bf C}^n \times {\bf C}^k \times {\bf C}^K$ с координатами $(z, \, w=u+i \, v, \, W=U+i \,V)$, заданная уравнениями
\begin{eqnarray}  \label{Q}
v = \Phi(z, \, \bar{z}), \quad  V=2 \, {\rm Re} \Psi(z,z,\bar{z}),
\end{eqnarray}
где $\Phi$ и $\Psi$ -- вектор-значные формы, линейные по каждому своему аргументу. В этой ситуации конечность типа равносильна тому, что координаты этих форм линейно независимы. Введем веса переменных следующим образом:
$[z]=1, \quad [w]=2, \quad [W]=3$.  Эта градуировка естественным образом распространяется на комплексные и вещественные степенные ряды. А с помощью дополнительного соглашения
$$[\frac{\partial}{\partial \, z}]=-1, \quad  [\frac{\partial}{\partial \, w}]=-2, \quad [\frac{\partial}{\partial \, W}]=-3$$
 и на векторные поля с аналитическими коэффициентами. После чего алгебра Ли инфинитезимальных голоморфных автоморфизмов в окрестности начала координат ${\rm aut}\, Q$ становится градуированной алгеброй Ли.  А если $Q$ -- невырождена, то конечномерной и конечноградуированной алгеброй вида
$${\rm aut} \, Q  =  g_{-3}+g_{-2}+ g_{-1}+g_{0}+g_{1}+\dots+g_{\delta}. $$
Эта алгебра образована векторными полями вида
$$X=2 \, {\rm Re} \left(f(z,w,W) \,\frac{\partial}{\partial z}+g(z,w,W) \,\frac{\partial}{\partial w}+h(z,w,W) \,\frac{\partial}{\partial W}\right),$$
где коэффициенты полей голоморфны в окрестности начала координат,  а поля удовлетворяют условию касания, т.е.
\begin{eqnarray}   \label{TQ}
{\rm Im} \, g=2  {\rm Re}  \Phi(f, \, \bar{z}), \quad  {\rm Im}\, h =2  {\rm Re} (2 \, \Psi(f,z,\bar{z}) +\Psi(z,z,\bar{f})),
\end{eqnarray}
при $w=u+i \,\Phi(z, \, \bar{z}), \; W=U+2 \,i \,{\rm Re} ( \Psi(z,z,\bar{z})$.

То что все $g_j$ при  $j>\delta$ равны нулю означает, в частности, что коэффициенты полей -- это полиномы, чьи степени не превышают $d=\delta+3$.

Ясно, что поля веса (-2) и (-3) - это поля вида
$$   X_{-3} =  2 \, {\rm Re} ( \mu \,\frac{\partial}{\partial W}),  \quad  X_{-2} =  2 \, {\rm Re} ( \nu \,\frac{\partial}{\partial w}),$$
где $\mu$ и $\nu$ -- произвольные постоянные вектора из ${\bf R}^K$ и ${\bf R}^k$.  Поле веса (-1) - это поле с коэффициентами весов ноль, один и два соответственно, т.е.
$$ f=const = p \in {\bf C}^n,    \quad  g=a (z,p),  \quad h=\alpha (z,z,p)+\beta (w,p). $$
Подставляя эти коэффициенты в условие касания, получаем
$$f = p,   \quad  g=2 \, i \,   \Phi(z, \, \bar{p}), \quad  h=2 \, i \,   \Psi(z,z, \, \bar{p}) + \beta(w,p),$$
где $\beta (w,p)$ -- это вещественная линейная форма, которая определяется из соотношения
\begin{equation}\label{H}
\beta (\Phi(z,\bar{z}),p)=4 \, {\rm Re} \Psi(p,z,\bar{z}).
\end{equation}
Единственность решения этого уравнения относительно $\beta$ гарантируется линейной независимостью коэффициентов
$\Phi$, а разрешимость при любом $p$ -- это условие, которое равносильно голоморфной однородности $Q$.\\

Подалгебре $g_{-}=g_{-3}+g_{-2}+g_{-1}$ соответствует подгруппа Ли $\mathcal{G}_{-}$, состоящая из треугольно-квадратичных сдвигов.   Голоморфная однородность $Q$ равносильна тому, что $\mathcal{G}_{-}$ действует на  $Q$ транзитивно, что позволяет отождествить  модельную поверхность $Q$ и группу Ли $\mathcal{G}_{-}$  как $CR$-многообразия. Заметим также, что при любых $\Phi$ и $\Psi$ компонента $g_0$ содержит поле вида
$$X_0=2 \, {\rm Re} \left(z \,\frac{\partial}{\partial z} + 2 \, w \,\frac{\partial}{\partial w}+3 \,W \,\frac{\partial}{\partial W}\right),$$

\vspace{3ex}

Пусть
$$\chi=(z \rightarrow F(z,w,W), \; w \rightarrow  G(z,w,W), \; W \rightarrow H(z,w,W))$$
 -- автоморфизм $Q$.  Заменив
$\chi$ на его композицию с соответствующим преобразованием из $\mathcal{G}_{-}$, можем считать, что  $\chi$ оставляет
начало координат неподвижным. Здесь мы используем голоморфную однородность $Q$.  Дифференциал автоморфизма $\chi$ переводит векторные поля в окрестности начала координат из ${\rm aut} \,  Q$ в векторные поля из ${\rm aut} \,  Q$.  Записывая, что векторные поля с координатами  $(A,B,C)$ и $(R,S,T)$ связаны отображением $\chi$, получаем
\begin{eqnarray}  \label{17}
\left[ \begin {array}{ccc} F_{{z}}&F_{{w}}&F_{{W}}\\
\noalign{\medskip}G_{{z}}&G_{{w}}&G_{{W}}\\
\noalign{\medskip}H_{{z}}&H_{{w}}&H_{{W}}\end {array} \right]^{-1}
 \cdot \left[ \begin {array}{c} A(F,G,H)
\\ \noalign{\medskip}B(F,G,H)\\ \noalign{\medskip}C(F,G,H)\end {array} \right]
=
 \left[ \begin {array}{c} R \\ \noalign{\medskip}S\\ \noalign{\medskip}T \end {array} \right]
\end{eqnarray}
Пусть $(e_1, \dots,e_n)$ -- стандартный базис пространства ${\bf C}^n$,  $(\nu_1, \dots,\nu_k)$ -- пространства ${\bf R}^k$,  $(\mu_1, \dots,\mu_K)$ -- пространства ${\bf R}^K$.   Выбирая эти значения для параметров, определяющих поля из $g_{-1},\,g_{-2},\,g_{-3}$ мы получаем поля, порождающие $g_{-}$ . Подставим все эти поля в (\ref{17}) вместо $(A,B,C)$,
а полученный результат запишем в блочно-матричной форме. Имеем
\begin{eqnarray}  \label{18}
\left[ \begin {array}{ccc} F_{{z}}&F_{{w}}&F_{{W}}\\
\noalign{\medskip}G_{{z}}&G_{{w}}&G_{{W}}\\
\noalign{\medskip}H_{{z}}&H_{{w}}&H_{{W}}\end {array} \right]^{-1}
 \cdot
 \left[ \begin {array}{ccc} E_{{n}}&{{0}}&{{0}}\\
\noalign{\medskip}2 i \Phi(F,E_n)&E_{{k}}&{{0}}\\
\noalign{\medskip}2i \Psi(F,F,E_n)+\beta(G,E_n)&{{0}}&E_{{K}}\end {array} \right]=P
\end{eqnarray}
где $(E_n,E_k,E_K)$ -- единичные матрицы соответствующих размеров,  $P$ -- матрица размера $ N \times N$,  составленная из векторов $(R,S,T)$ -- образов векторов из  $g_{-}$, т.е. её элементы -- это полиномы, степени не выше $d$.
В силу соотношения (\ref{18}) эта матрица -- невырождена и элементы
обратной матрицы $M=P^{-1}$ -- это рациональные функции, степени не выше $d \, N$. Пусть
\begin{eqnarray*}
M=P^{-1}=
\left[ \begin {array}{ccc} M_{{1}}^1&M_{{2}}^1&M_{{3}}^1\\
\noalign{\medskip}M_{{1}}^2&M_{{2}}^2&M_{{3}}^2\\
\noalign{\medskip}M_{{1}}^3&M_{{2}}^3&M_{{3}}^3\end {array} \right].
\end{eqnarray*}
Запишем (\ref{18}) в виде
\begin{eqnarray}  \label{19}
\left[ \begin {array}{ccc} F_{{z}}&F_{{w}}&F_{{W}}\\
\noalign{\medskip}G_{{z}}&G_{{w}}&G_{{W}}\\
\noalign{\medskip}H_{{z}}&H_{{w}}&H_{{W}}\end {array} \right]
 =\\
\nonumber
 \left[ \begin {array}{ccc} E_{{n}}&{{0}}&{{0}}\\
\noalign{\medskip}2 i \Phi(F,E_n)&E_{{k}}&{{0}}\\
\noalign{\medskip}2i \Psi(F,F,E_n)+\beta(G,E_n)&{{0}}&E_{{K}}\end {array} \right] \cdot
\left[ \begin {array}{ccc} M_{{1}}^1&M_{{2}}^1&M_{{3}}^1\\
\noalign{\medskip}M_{{1}}^2&M_{{2}}^2&M_{{3}}^2\\
\noalign{\medskip}M_{{1}}^3&M_{{2}}^3&M_{{3}}^3\end {array} \right]
\end{eqnarray}
Из первой блок-строки этого соотношения получаем
$${\rm grad}\, F=(F_z,F_w,F_W)=(M_1^1,M_2^1,M_3^1),$$
т.е. ${\rm grad}\, F$ -- рационален, степени не выше $d \,N$.

Степень рациональной функции -- это максимум степеней числителя и знаменателя.  Аналогично и с весом (взвешенной степенью). Для дальнейших подсчетов отметим, что любые арифметические операции с двумя рациональными функциями
степеней $d_1$ и $d_2$ дают рациональную функцию, чья степень не превосходит $d_1+d_2$.

\vspace{2ex}

{\bf Лемма 1:}  Пусть ${\rm deg} \, R_1 =d_1, \; {\rm deg} \, R_2 =d_2$, тогда ${\rm deg} \, (R_1 \diamond R_2) \leq (d_1+d_2)$, где $\diamond$ -- любая из четырех арифметических операций. То же самое верно и для веса рациональной функции.

\vspace{2ex}

Подставим в  (\ref{18}) в качестве $(A,B,C)$ поле $X_0$,
получаем
\begin{eqnarray}  \label{20}
 \left[ \begin {array}{c} F
\\ \noalign{\medskip}2 \, G\\ \noalign{\medskip}3 \, H\end {array} \right]=
\left[ \begin {array}{ccc} F_{{z}}&F_{{w}}&F_{{W}}\\
\noalign{\medskip}G_{{z}}&G_{{w}}&G_{{W}}\\
\noalign{\medskip}H_{{z}}&H_{{w}}&H_{{W}}\end {array} \right]
\left[ \begin {array}{c} R \\ \noalign{\medskip}S\\ \noalign{\medskip}T \end {array} \right]
\end{eqnarray}
Первая блок-координата этого соотношения имеет вид  $F=F_z \, R + F_w \, S +F_W \, T$, откуда получаем, что $F$ -- рациональна, степени не выше $d \,(N+1)$.
Возвращаясь к (\ref{19}), из второй блок-строки получаем
\begin{eqnarray*}
G_z= 2 i \Phi(F,E_n) M_1^1+M_1^2,\\
G_w= 2 i \Phi(F,E_n) M_2^1+M_2^2,\\
G_W= 2 i \Phi(F,E_n) M_3^1+M_3^2.\\
\end{eqnarray*}
Следовательно  ${\rm grad}\, G$ -- рационален, степени не выше $d \,(3N+1)$.  Теперь из второй координаты (\ref{20}) следует, что $G$ -- рациональна, степени не выше $d \,(3 N+2)$. Аналогично из третьей блок строки (\ref{19}) получаем
\begin{eqnarray*}
H_z= (2 i \Psi(F,F,E_n)+\beta(G,E_n)) M_1^1+M_1^3,\\
H_w=(2 i \Psi(F,F,E_n)+\beta(G,E_n)) M_2^1+M_2^3,\\
H_W= (2 i \Psi(F,F,E_n)+\beta(G,E_n)) M_3^1+M_3^3.\\
\end{eqnarray*}
Откуда следует, что ${\rm grad}\, H$ -- рационален, степени не выше $d \,(7N+6)$.  Теперь из третьей координаты (\ref{20}) следует, что $H$ -- рациональна, степени не выше $7 \, d \,( N+1)$.

Итак, любой автоморфизм $\chi$, сохраняющий начало координат -- это рациональное отображение, степени не выше $7 \, d \,( N+1)$.  Произвольный автоморфизм $Q$ имеет вид $\eta(\chi)$, где $\eta \in \mathcal{G}_{-}$ -- квадратично-треугольное преобразование. Поэтому итоговая оценка степени для произвольного автоморфизма --  $14 \, d \,( N+1)$.

\vspace{2ex}

{\bf Лемма 2:} Пусть $Q \in \mathbf{C}^N$ -- невырожденная модельная поверхность.
Тогда коэффициенты векторных полей, составляющих ${\rm aut}\,Q$ -- это полиномы степени не выше $N^3$. \\
{\it Доказательство:} Пусть $Q $ -- модельная поверхность
$CR$-размерности $n$  и коразмерности $\kappa$ ($N=n+\kappa$) и  ${\rm aut}\,Q=g_{-l}+ \dots +g_0 +\dots +g_{\delta}$, т.е. $\delta$ - это старший вес.
Из результатов работы \cite{Z} непосредственно следует, что  автоморфизм  $Q$ однозначно определяется своей $n \,(\kappa+1)$-струей в точке. Отсюда следует, что то же самое можно утверждать и о коэффициентах полей из ${\rm aut}\,Q$. Поскольку вместе с любым векторным полем $X$ алгебра ${\rm aut}\,Q$ содержит каждую его градуированную компоненту, то она не может содержать весовых компонент, чей вес превосходит $l \,  n \,(\kappa+1)$,
т.е. $d \leq \delta \leq l \,  n \,(\kappa+1) \leq   n \,\kappa  \,(\kappa+1) \leq N^3$.

\vspace{2ex}
В частности, наше рассуждение показывает, что для поверхности $Q$ типа $m=((2,k),(3,K))$ имеем $d \leq 3 \, n \, (k+K+1)$.
\vspace{2ex}

{\bf Утверждение 3 :} Пусть $Q$ -- невырожденная голоморфно однородная модельная поверхность $CR$-размерности $n$ Блум-Грэм-типа $m=((2,k),(3,K))$,  коразмерности $k+K$, тогда ${\rm Aut} \, Q$ состоит из бирациональных преобразований пространства $\mathbf{C}^N$,   где $N=n+k+K$, чья степень  не превосходит
$$D(n,k,K) \leq 42 \, n \,(k+K+1)(n+k+K+1)   \;\;  \mbox{или} \;\;     D(N) \leq  \frac{21}{2} \,(N^2-1)(N+3).$$

\vspace{5ex}

Пусть теперь $Q$ -- это произвольная невырожденная голоморфно однородная модельная поверхность. Как было показано в \cite{VB20} набор весов при условии однородности -- это отрезок натурального ряда $(2,3,\dots,l)$.  В таком случае координаты объемлющего пространства нам будет удобно обозначить
$$ (z,w_2,w_3,\dots, w_l), \quad z \in \mathbf{C}^n, \; w_j=u_j+i\,v_j \in \mathbf{C}^{k_j}.$$
При этом веса переменным назначены так: $[z]=[\bar{z}]=1, \, [w_j]=[u_j]=j, \; j=2,\dots,l$. Это соглашение вводит градуировку степенных рядов и векторных полей. Теперь уравнения $Q$ можно записать так:
\begin{eqnarray} \label{Q}
  v_j = \Phi_j(z,\bar{z},u_2,\dots,u_{j-1}), \quad j=2,\dots,l
\end{eqnarray}
где вещественная вектор-значная форма $\Phi_j$ однородна веса $j$ и записана в приведенной форме (см. \cite{VB20}).
В силу невырожденности алгебра ${\rm aut} \, Q =g_{-}+g_0+g_{+}$ конечномерна, конечноградуирована и состоит из полей с полиномиальными коэффициентами равномерно ограниченной степени.  Условие голоморфной однородности $Q$ равносильно тому, что ${\rm  dim} \,  g_{-} = {\rm  dim} \, Q$. Рассмотрим подробнее строение  подалгебры
$ g_{-}=g_{-l}+g_{-l+1}+\dots+g_{-2}+g_{-1}$. Она состоит из полей вида
$$X=2 \, {\rm Re} \left(f(z,w_2,\dots,w_l ) \,\frac{\partial}{\partial z}+g_2(z,w_2,\dots,w_l) \,\frac{\partial}{\partial w_2
}+\dots + g_l(z,w_2,\dots,w_l) \,\frac{\partial}{\partial w_l}\right),$$
удовлетворяющих условию касания, а именно
\begin{eqnarray*}
{\rm Im} \,g_j = d \Phi_j (z,\bar{z},u_2,\dots,u_{j-1})(f,\bar{f}, {\rm Re}( \,g_2),\dots, {\rm Re}( \,g_l)),\\
\mbox{при  } w_{j}=u_{j}+i \, \Phi_{j}( z,\bar{z},u_2,\dots,u_{j-1}), \; j=2,\dots,l.
\end{eqnarray*}
Или просто в координатах $X=(f(z,w),\; g_2(z,w), \dots, g_l(x,w))$.  Если $X_j \in g_{j}$, то пишем
$$X_j=(f_j(z,w),g_{2,j}(z,w),\dots,g_{l,j}(z,w)).$$
При этом $[f_j]=j+1, \; [g_{\nu,j}]=\nu+j$, а степени, соответственно, не превосходят весов.  Для всех $-l \leq j \leq -2$
\begin{eqnarray*}
f_j(z,w)=g_{j,2}(z,w)= \dots = g_{j,-j-1}(z,w)=0, \\
g_{j,-j}(z,w)=\beta_j \in \mathbf{R}^{k_j}, \quad
 g_{j,\nu}(z,w)=g_{j,\nu}(z,w,\beta_j),   \;  - j +1  \leq \nu \leq l
\end{eqnarray*}
Где  $g_{j,\nu}(z,w,\beta_j)$ -- это полином по $(z,w)$ веса $j+\nu$,  линейный по $\beta_j$.
Соответственно для $j=-1$ имеем
\begin{eqnarray*}
f_{-1} =p \in \mathbf{C}^{n},  \quad  g_{-1,\nu}(z,w)=g_{-1,\nu}(z,w,p),   \quad  2  \leq \nu \leq l
\end{eqnarray*}
Где  $g_{-1,\nu}(z,w,p)$ -- это полином по $(z,w)$ веса $\nu-1$,  вещественно линейный по $p$.
То что условия
касания  однозначно разрешимы относительно $g_{j,\nu}(z,w,\beta_j)$ при любом фиксированном $\beta_{j}$ и
относительно  $g_{-1,\nu}(z,w,p)$ при фиксированном $p$ -- это прямое следствие голоморфной однородности $Q$.
Заметим также, что при любых $\Phi_j$ компонента $g_0$ содержит поле вида
\begin{equation}\label{g0}
X_0=2 \, {\rm Re} \left(z \,\frac{\partial}{\partial z} + 2 \, w_2 \,\frac{\partial}{\partial w_2}+\dots + l \,w_l \,\frac{\partial}{\partial w_l}\right)
\end{equation}

\vspace{3ex}
Переходя к доказательству бирациональности, отметим, что мы следуем той же схеме Каупа-Туманова, которая была продемонстрирована выше.   Пусть
$$\chi=(z \rightarrow F(z,w_2,\dots,w_l), \quad  w_j  \rightarrow  G^j (z,w_2,\dots,w_l)), \quad j=2,\dots,l $$
 -- автоморфизм $Q$.  Заменив $\chi$ на его композицию с соответствующим преобразованием из $\mathcal{G}_{-}$, можем считать, что  $\chi$ оставляет начало координат неподвижным. Дифференциал автоморфизма $\chi$ переводит любые векторные поля в окрестности начала координат из ${\rm aut} \,  Q$ в векторные поля из ${\rm aut} \,  Q$ .
 Записывая, что векторные поля с координатами  $(A,B_2,\dots,B_l)$ и $(R,S_2,\dots,S_l)$ связаны отображением $\chi$, получаем
\begin{eqnarray}  \label{27}
\left[ \begin {array}{cccc} F_{{z}}&F_{{w_2}}& \dots & F_{{w_l}}\\
\noalign{\medskip} G^2_{{z}}&G^2_{{w_2}}& \dots & G^2_{{w_l}}\\
\noalign{\medskip} \dots&\dots& \dots & \dots\\
\noalign{\medskip} G^l_{{z}}&G^l_{{w_2}}& \dots & G^l_{{w_l}}\\
\end {array} \right]^{-1}
 \cdot \left[ \begin {array}{c} A(F,G^2,\dots,G^l)
\\ \noalign{\medskip}B_2(F,G^2,\dots,G^l)\\
\\ \noalign{\medskip}\dots\\
\\ \noalign{\medskip}B_l(F,G^2,\dots,G^l)\\  \end {array} \right]
=
 \left[ \begin {array}{c} R \\ \noalign{\medskip}S_2\\ \noalign{\medskip} \dots \\
 \noalign{\medskip}S_l \end {array} \right]
\end{eqnarray}
При этом, как и раньше, мы пользуемся блочно-матричной арифметикой, т.е. квадратную матрицу размером $N \times N$ мы записываем как блочную $l \times l$ матрицу.

Пусть $e=(e_1, \dots,e_n)$ -- стандартный базис пространства ${\bf C}^n$ и пусть   $\nu^j=(\nu^j_1, \dots,\nu^j_{k_j})$ --  базис пространства ${\bf R}^{k_j}$.   Выбирая элементы $e$ в качестве значений для параметров, опредедяющих поля из $g_{-1}$ а $\nu^j$ -- поля из $g_{-j}$, мы получаем поля, порождающие всю подалгебру $g_{-}$ . Располагая все такие поля в виде столбцов блочной квадратной матрицы $T$, получаем
\begin{eqnarray*}
\left[ \begin {array}{ccccc} E_n & 0 & 0 &  \dots & 0\\
\noalign{\medskip} g_{-1,2}(F,E_n)&E_{k_2}& 0& \dots & 0\\
\noalign{\medskip} \dots&\dots& \dots & \dots& 0\\
\noalign{\medskip} g_{-1,l-1}(F,G^2, \dots,E_n)&g_{-2,l-1}(F,G^2, \dots,E_{k_2}) & \dots & E_{k_{l-1}}& 0\\
\noalign{\medskip} g_{-1,l}(F,G^2, \dots,G^{l-1},E_n)&g_{-2,l}(F,G^2, \dots,E_{k_2})&g_{-3,l}(F,G^2, \dots,E_{k_3})& \dots & E_{k_l}\\
\end {array} \right]
 \end{eqnarray*}
 Обозначим блок-элемент матрицы $M$, стоящий на пресечении $i$-й блок-строки и $j$-го блок столбца через $M^i_j$, т.е.
\begin{eqnarray*}
M=P^{-1}=
\left[ \begin {array}{cccc} M_{{1}}^1&M_{{2}}^1&\dots&M_{{l}}^1\\
\noalign{\medskip}M_{{1}}^2&M_{{2}}^2&\dots &M_{{l}}^2\\
\noalign{\medskip}\dots&\dots&\dots &\dots\\
\noalign{\medskip}M_{{1}}^l&M_{{2}}^l&\dots&M_{{l}}^l\end {array} \right].
\end{eqnarray*}
Теперь полученные из (\ref{27}) соотношения запишем в виде одного матричного равенства
\begin{equation}\label{28}
 J = T \cdot P^{-1}=T \cdot M
\end{equation}
Где $J$ -- якобиева матрица отображения $\chi$,  а $P$ -- матрица, составленная из полей вида  $(R,S_2,\dots, S_l)$ -- образов базисных полей из $g_{-1}$.  Если число $d$ оценивает сверху степени полей из ${\rm aut} \, Q$, то элементы матрицы $M=P^{-1}$ --  это рациональные функции, чьи степени не превосходят $N \, d$.

Записывая первую блок-строку соотношения (\ref{28}), получаем
\begin{equation}\label{29}
{\rm grad}\, F=(F_z,F_{w_2}, \dots, F_{w_l})=(M_1^1,M_2^1,\dots,M_l^1),
\end{equation}
т.е. ${\rm grad}\, F$ -- рационален, степени не выше $d \,N$.
Подставим, далее,  в  (\ref{27}) в качестве $(A,B_2,\dots,B_l)$ поле $X_0 \in g_0$ (см.(\ref{g0})),
получаем
\begin{eqnarray}  \label{30}
 \left[ \begin {array}{c} F
\\ \noalign{\medskip}2 \, G^2\\ \noalign{\medskip}\dots \\ \noalign{\medskip}l \, G^l \end {array} \right]=
\left[ \begin {array}{cccc} F_{{z}}&F_{{w_2}}&\dots &F_{{w_l}}\\
\noalign{\medskip}G^2_{{z}}&G^2_{{w_2}}&\dots &G^2_{{w_l}}\\
\noalign{\medskip}\dots&\dots&\dots &\dots\\
\noalign{\medskip}G^l_{{z}}&G^l_{{w_2}}&\dots&G^l_{{w_l}}\end {array} \right]
\left[ \begin {array}{c} R \\ \noalign{\medskip}S_2\\ \noalign{\medskip}\dots \\ \noalign{\medskip}S_l \end {array} \right]
\end{eqnarray}
Записывая первую блок-координату (\ref{30}), из (\ref{29}) получаем $F=F_z \, R + F_{w_2} \, S_2 +\dots +F_{w_l} \, S_l$.
 Поскольку $l \leq K \leq N$, то  мы заключаем, что ${\rm deg} \, F \leq 2 \, d \, N^2$.

Теперь, записывая вторую блок-строку  (\ref{28}), получаем
\begin{eqnarray*}
G^2_z= g_{-1,2} (F,E_n) M_1^1+M_1^2,\\
G^2_{w_2}= g_{-1,2}(F,E_n) M_2^1+M_2^2,\\
\dots\\
G^2_{w_l}= g_{-1,2}(F,E_n) M_l^1+M_l^2.\\
\end{eqnarray*}
Следовательно  ${\rm grad}\, G^2$ -- рационален, степени не выше $d \,( 2 \, N^2+2 \, N) $.
Тогда из второй координаты (\ref{30}) следует, что $G^2$ -- рациональна, степени не выше $d \,( 2 \, N^2+2 \, N+1)\,N$.

Оценка степени $G^2$  приведена здесь для удобства читателя. Она включается в общее рекуррентное рассуждение, приведенное ниже.

Итак, теперь мы готовы описать общий ($j+1$)-й этап этого процесса -- оценку степени $G^{j+1}$. При этом мы будем давать верхнюю оценку степени, зависящую только от $N \geq 2$, причем не претендующую на точность. Этот этап, как мы видели состоит из двух шагов:  оценка степени градиента $G^{j+1}$, основанная на соотношении   (\ref{28}),  и оценка степени самой компоненты $G^{j+1}$,  основанная на соотношении   (\ref{30}).

Итак, пусть $d_1$ -- это полученная нами оценка степени $F$, т.е. $d_1 = d  \, \sigma_1 =2 \, d \, N^2$ и пусть на $j$-м
этапе получена величина $d_j=d \, \sigma_j$ -- оценка степени $G^j$. Рассмотрим $(j+1)$-ю блок-строку соотношения  (\ref{28}).  В левой части такого соотношения стоит ${\rm grad}\, G^{j+1}$, в правой его части -- выражение в виде суммы, где число слагаемых не превосходит $l \leq K \leq N$. При этом каждое слагаемое имеет вид произведения
$g_{\nu,j+1}(F,G^2,\dots,E_{k_{-\nu}})$ и некоторого блок-элемента матрицы $M$. Учитывая, что степени  $ g_{\nu,j+1}$ меньше $l$,  а степени элементов $M$ не выше $d \,N$,  мы можем утверждать, что степень ${\rm grad}\, G^{j+1}$ не превосходит $d \, (l \, \sigma_j + N)\,l$.   Далее, выписывая $(j+1)$-ю блок-координату  (\ref{30}), получаем, что $G^{j+1}$ -- это сумма длины $l$, где степени слагаемых не превосходят
$ d \, (l \, \sigma_j + N)\,l +d$, откуда следует, что $d \, \sigma_{j+1} \leq   d \, l\,( (l \, \sigma_j + N)\,l +1).$
Таким образом
$$\sigma_{j+1} \leq l \, ((l \, \sigma_j +N)\,l+1) \leq N^3 \, (\sigma_j+2) \leq 2 \, N^3 \, \sigma_j.$$
Учитывая, что $\sigma_1 = 2 \, N^2 \leq 2 \, N^3$, и то, что длина последовательности $(\sigma_1,\sigma_2, \dots)$ не превосходит $N$, получаем, что все степени не превосходят $d \, (2\, N^3)^N$.  В соответствии с леммой 2 имеем
$d  \leq N^3$.   Поэтому величина $N^3 \, (2\, N^3)^N$ даёт оценку для степени автоморфизма $\chi$.
Теперь следует вспомнить, что произвольный автоморфизм $Q$ можно получить из $\chi$ композицией с треугольно полиномиальным преобразованием из $\mathcal{G}_{-}$, чья степень не превосходит $l -1\leq N$. В итоге получаем

\vspace{2ex}

{\bf Теорема 4:}  Группа ${\rm Aut} \, Q$ голоморфных автоморфизмов произвольной невырожденной голоморфно однородной модельной поверхности $Q \subset \mathbf{C}^N$ состоит из биголоморфных преобразований $\mathbf{C}^N$, чья максимальная степень не превосходит
$$ D(N) \leq N^4 \, 2^N \, (N^N)^3.$$

\vspace{2ex}
В теореме 4 имеется два условия: невырожденность и голоморфная однородность.  При этом условие невырожденности, в свою очередь, распадается на два условия: конечность типа и голоморфную невырожденность.
Если $Q$ - это модельная поверхность, чей тип по Блуму-Грэму -- бесконечен, то она может быть задана в виде
(\ref{Q}),  где среди координат какой-либо взвешенно-однородной формы $\Phi_j$ стоит тождественный ноль.
Тогда подвергая соответствующую координату  группы $w$ произвольному вещественно-аналитическому преобразованию (остальные координаты -- на месте), получаем автоморфизм $Q$. Т.е. в этом случае утверждение о бирациональности -- неверно.  Оно также становится неверным, если нарушено второе условие невырожденности - голоморфная невырожденность. Т.к. в этом случае в алгебре  ${\rm aut} \, Q$ имеются не полиномиальные поля,
что противоречит бирациональности ${\rm Aut} \, Q$.

Покажем, что голоморфная однородность также является необходимым условием.  Соответствующий пример был рассмотрен в \cite{VB97}.

\vspace{2ex}

{\bf Пример 5:}  Рассмотрим модельную гиперповерхность пространства $\mathbf{C}^2$ вида
$Q=\{v=|z|^4\}$.   Эта гиперповерхность не является голоморфно однородной. Её Блум-Грэм-тип меняется от точки к точке.  Пусть $\xi=(a,b) \in Q$ и $m(\xi)$ - тип в точке $\xi$,  тогда
\begin{eqnarray*}
m((a,b))= (4),  \mbox{  если }  \; a =0  \\
m((a,b))= (2),  \mbox{  если }  \; a \neq 0  \\
\end{eqnarray*}
При этом ${\rm aut} \, Q$ содержит поле
$$X={\rm Re} \left( w \, z \,\frac{\partial}{\partial \, z}+ 2 \, w^2 \,  \frac{\partial}{\partial \, w}\right).$$
Этому полю соответствует 1-параметрическая подгруппа ${\rm Aut} \, Q$ преобразований вида
$$ z \rightarrow \frac{z}{\sqrt{1-2\,t\, w}},    \quad    w \rightarrow \frac{w}{1-2\,t\, w},$$
которые не являются рациональными.

\vspace{2ex}

Модельные поверхности, очевидно, попадают в класс вещественно алгебраических многообразий.  Поэтому можно утверждать, что автоморфизмы любой невырожденной модельной поверхности, не зависимо от голоморфной однородности, являются {\it алгебраическими} (см. \cite{BER}, теорема 13.1.4.).

\vspace{2ex}

Группы автоморфизмов невырожденных однородных модельных поверхностей демонстрируют примеры подгрупп группы бирациональных автоморфизмов комплексного аффинного пространства с условием равномерной ограниченности степеней.  Такие группы представляют интерес не зависимо от $CR$-геометрии. В работе
\cite{ZH} приводится красивая конструкция позволяющая строить такие группы (см.теорему 3). К сожалению, она не применима в нашей ситуации.  В связи с этим возникает запрос на некую, более общую конструкцию, которая могла бы включить и группы автоморфизмов модельных поверхностей.

\vspace{3ex}
К вопросу об оценке степени можно подойти более дифференцированно. Например, можно рассмотреть величину $D(m)$, которая представляет собой максимум степеней автоморфизмов по всем модельным поверхностям фиксированного конечного Блум-Грэм-типа $m$.
Из теоремы 4  сразу следует, что эта величина конечна.  Далее, можно, по крайней мере, в простейших ситуациях, ставить вопрос о точных значениях $D(m)$.    Например, если речь идет о модельных гиперквадриках ($m=(2)$),   то $D((2))=1$.   Если же это квадратичные модельные поверхности более высокой коразмерности $k>1$ (квадрики),  т.е. $m=(2,...,2)=(2,k)$, то $D((2,k)) \geq k$.   Примеры квадрик с автоморфизмами, чья степень больше коразмерности, нам не известны.     

\vspace{5ex}

{\bf  О топологическом строении группы ${\rm Aut} \, Q_{0}$}.  Отметим сразу, что в отличие от исследования вопроса о бирациональности в предыдущем разделе, в этом разделе  голоморфная однородность $Q$ предполагаться не будет.

Как было показано в \cite{VB20},  если $M_{\xi}$ - росток невырожденной {\it вещественно алгебраической} поверхности,  то ${\rm Aut} \, M_{\xi}$ состоит из алгебраических отображений, голоморфных в окрестности $\xi$  (см. \cite{ZH}) и обладает структурой группы Ли (\cite{BER}, теорема 12.7.18).   В частности, этот результат применим к любой невырожденной модельной поверхности $Q$.

Итак, пусть $Q$ - невырожденная модельная поверхность. Введем в рассмотрение следующие объекты:\\
$g_{-}$ - подалгебра  ${\rm aut} \, Q$ состоящая из полей отрицательного веса,  $G_{-}$ - соответствующая $g_{-}$  связная группа Ли;\\
$g_{0}$ - подалгебра  ${\rm aut} \, Q$, состоящая из полей веса ноль ,  $G_{0}$ - соответствующая $g_{0}$  связная группа Ли;\\
$g_{+}$ - подалгебра  ${\rm aut} \, Q$, состоящая из полей положительного веса ,  $G_{+}$ - соответствующая $g_{+}$  связная группа Ли;\\
$St$ - стабилизатор начала координат в ${\rm Aut} \, Q_{0}$

$\mathcal{G}_{-}$ - это подгруппа треугольно полиномиальных автоморфизмов $Q$, описанная в \cite{VB20}.\\
$\mathcal{G}_{0}$ - это подгруппа автоморфизмов $Q$,  описанная в \cite{VB20}, т.ч.  действие на координату $z$ имеет вид $(z \rightarrow C \, z)$, где $C$ - невырожденное линейное преобразование. \\
$\mathcal{G}_{+}$ - это подгруппа автоморфизмов $Q$ вида  $(z \rightarrow  z + o(1), \; w_j \rightarrow w_j +o(m_j)), \;  j=1,...,l.$

\vspace{3ex}

{\bf Утверждение 6:}  (a)  Любой автоморфизм $Q$ представим в виде $\tau \circ \sigma$, где $\sigma \in St, \;  \tau \in \mathcal{G}_{-} $. При этом $\sigma$ и $\tau$ определены однозначно и имеет место полупрямое разложение  ${\rm Aut} \, Q=  \mathcal{G}_{-} \ltimes St$.\\
(b)   Любой автоморфизм $Q$  из $St$ представим в виде $L \circ N$, где $L \in \mathcal{G}_{0}, \;  N \in \mathcal{G}_{+} $. При этом $L$ и $N$ определены однозначно и имеет место полупрямое разложение  $St=  \mathcal{G}_{0} \ltimes \mathcal{G}_{+}$.\\
(c)   Группа $\mathcal{G}_{0}$ квазилинейных преобразований $(z \rightarrow  C \, z , \; w_j \rightarrow \rho_j(w)),  \;j=1,..., l,$  имеет точное представление вида $(z \rightarrow  C \, z , \; w_j \rightarrow \rho_j(w))  \rightarrow (z \rightarrow  C \, z )$  в $GL(n,{\bf C})$ и изоморфна
вещественной линейной алгебраической группе.  \\
{\it Доказательство:}  Разложения из (a) и (b) обсуждались в  \cite{VB20}. Утверждения о полупрямом произведении проверяются непосредственно. Точность представления из (c) - это  Теорема 5, п.(f) \cite{VB20}. Утверждение доказано.

  \vspace{3ex}

{\bf Теорема 7:} (a)  $\mathcal{G}_{+}=G_{+}$, в частности группа Ли $\mathcal{G}_{+}$ - связна и односвязна;\\
(b) Группа Ли $\mathcal{G}_{0}$ имеет конечное число компонент связности;\\
(с) $\mathcal{G}_{-}=G_{-}$, в частности группа Ли $\mathcal{G}_{-}$ - связна .\\
{\it Доказательство:}  Рекуррентный процесс вычисления компонент  $\mathcal{G}_{+}$, описанный  в \cite{VB20} (конструкция Пуанкаре) позволяет однозначно восстановить элемент $\mathcal{G}_{+}$ по параметрам, содержащимся в ядре гомологического оператора.  Совокупность этих параметров - это линейное пространство. Все параметры, принадлежащие ядру гомологического оператора,  реализуются элементами $G_{+}$. Пункт (а) - доказан.\\
В силу п.(с) утверждения 6 группа $\mathcal{G}_{0}$ изоморфна вещественной линейной алгебраической группе.  Пункт (b) - доказан.\\
Группа $\mathcal{G}_{-}$ состоит из преобразований $S_{\xi}$,  которые однозначно определяются выбором точки $\xi=(a,b_1, ...,b_l)$ того же Блум-Грэм-типа, что и начало координат.  Процесс построения $S_{\xi}$ описан в  \cite{VB20} (см. доказательство теоремы 17).  По существу это построение - это процесс приведения поверхности к стандартной форме  в точке $\xi$. В координатах связанных с  точкой $\xi$ каждая координатная форма имеет вид суммы весовых компонент. Старшая компонента совпадает с формой в начале координат, а младшие зависят от точки $\xi$. То, что точка $\xi$ является точкой того же типа, как и начало координат означает, что в процессе приведения все компоненты младших весов редуцируются к нулю. При этом из построения видно, что компонента каждого веса редуцируется отдельно. Поэтому преобразование
$a \rightarrow t \, a, \; b_j \rightarrow t^{m_j} \, b_j,  \; t>0, $ переводит $\xi$ в точку того же  Блум-Грэм-типа.  Устремляя $t$ к нулю получаем, что множество точек, имеющих тот же тип, что и начало координат - связно. Из этого следует стягиваемость и, в частности, связность $\mathcal{G}_{-}$.  Теорема доказана.

  \vspace{3ex}

{\bf Следствие 8:}  Группа ${\rm Aut} \, Q_{0}$ стягиваема на $\mathcal{G}_{0}$, поэтому число компонент связности и все гомотопические группы у них совпадают.

 \vspace{3ex}

{\bf Пример 9:}  Простейший пример модельной поверхности с несвязной $\mathcal{G}_{0}$ - это гиперквадрика в ${\bf C}^3$ со знаконеопределенной формой Леви  $\{ {\rm Im} w = |z_1|^2 -|z_2|^2\}$.  В этом примере  $\mathcal{G}_{0}$ имеет две компоненты.

\end{document}